\documentclass[final]{article}
\pdfoutput=1
\usepackage{graphicx}
\usepackage{amsmath,amssymb}
\usepackage{natbib}

\def\be{\begin{equation}}
\def\ee{\end{equation}}
\def\~{\mathaccent "7E}

\newtheorem{Lemma}{Lemma}[section]

\newtheorem{Remark}{Remark}[section]
\newtheorem{Proof}{Proof}[section]

\title{ A modified Allen-Cahn model for pattern synthesis on surfaces}
\author{Lorina Dascal \qquad 
Gautam Pai \qquad  Ron Kimmel \\[12pt]
{Technion - Israel Institute of Technology}\\[2pt]
\tt{\{lorina,paigautam,ron\}@cs.technion.ac.il}
}
\date{} 
\begin{document}
\maketitle
\begin{figure}[h!]
 \begin{center}
   \mbox{
   \includegraphics[width=0.45\textwidth,height=5cm] 
   {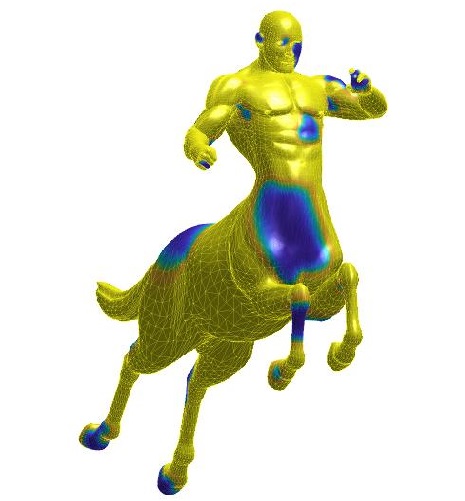}
   \includegraphics[width=0.45\textwidth,height=5cm] 
   {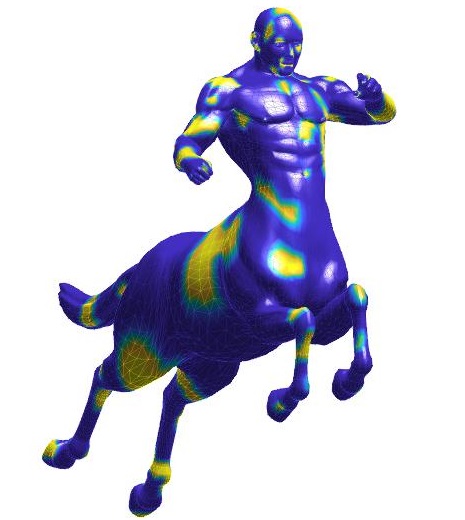}
   }
   \caption{
   Various patterns  on the surface by means of modified Allen-Cahn model. Left:Inverted spots. Right: Spots.
   }
   \end{center}
\end{figure}
\begin{abstract}
We propose an extension of the Allen-Cahn model for pattern synthesis on two dimensional curved surfaces. 
   This model is based on a single PDE and it offers improved  ability of controlling the type of generated surface patterns via the chosen  reaction-diffusion coefficient, thus, obtaining patterns in form of spots, inverted spots, or stripes. We investigate the   dependence of the type of the obtained pattern on the  new  proposed reaction term.
An efficient operator splitting scheme is used to discretize the  model on a surface. 
Experiments on surfaces with varying initial conditions illustrate  a variety of  patterns.
\end{abstract}
  \section{Introduction}
The   Allen-Cahn equation is a semilinear reaction diffusion based PDE, and  has been widely used in phase separation, crystal-growth  and various material science applications defined on planar domains.  
It was originally introduced in \cite{ACahn} as  a model   for the motion  of anti-phase  boundaries  in a binary alloy.  It was later used in phase separation \cite{phasesep}, or crystal growth analysis \cite{crystmod}. 
A study of  a different semi-linear reaction-diffusion equation, governed by the Hamiltonian operator was proposed in \cite{Hamilt}, altogether with applications for 3D spectral mesh geometry compression.
In this paper, we explore the behaviour of an extension of the Allen-Cahn reaction-diffusion partial differential equation on  triangulated surfaces and analyze  its generated patterns, which will be shown to exhibit a richer set of patterns compared to the regular Allen-Cahn equation itself.
Turing models and in general reaction-diffusion systems have been extensively used for understanding spatial patterns.  
Patterns arising in reaction-diffusion processes have been proposed in biology applications to describe developmental processes such as skin pigmentation patterning \cite{skin}.
The solvers used in various publications involve more advanced  numerical methods for solving the reaction-diffusion system, especially when applied on a curved domain. 
Allen-Cahn PDE was applied to two dimensional flat domain, as well as to curved surfaces, for which various numerical solvers were analyzed in \cite{RBF}, \cite{Dzk}. 
However, the only pattern that can be obtained with this model is a stripe-based structure.  
Its  complement in the large class of reaction diffusion systems  is the Fitzhugh-Nagumo reaction-diffusion system, see \cite{Nagumo}, which contains in its nonlinear part a third order polynomial.
It can generate patterns like  spots/stripes, but at the cost of solving a nonlinear coupled PDE system. 
In this paper, we propose a simple extension of the Allen-Cahn model, exhibiting richer patterns on general non-flat geometries by means of a single PDE defined on the surface.
We analyze the behavior of the suggested model with its newly introduced reaction term and handle it numerically by an efficient  operator-splitting scheme.  
Numerical examples illustrate the generation of spot and stripe patterns on surfaces as a function of the reaction parameter.   
Understanding the nature of the patterns obtained on various curved geometries is challenging. 
We purpose to analyze the relationship between the given coefficients of the governing reaction-diffusion equation, the underlying geometry, and the type of the resulting pattern. 
Moreover, the simulated numerical solution of the modified model is compactly supported, unlike the original Allen-Cahn model. Numerical experiments illustrate this important locality property.

\section{Modified Allen-Cahn equation on  surfaces}
 We describe below the modified Allen-Cahn  on surfaces. 
 Suppose  $\Gamma$ is a surface in $\mathbb{R}^3$, and $\partial \Gamma $ is empty.  
 The  known surface Allen-Cahn (AC) equation is as follows:
\begin{equation}
u_t =\Delta_{\Gamma} u-\frac{1}{\epsilon^2}f(u), \,\,x\in \Gamma, \,\, t\in[0,T]
\label{AC1}
\end{equation}
\begin{equation}
u|_{t=0}=u_0(x), x\in \Gamma
\end{equation}
where $ f(u)=u^3-u$,  $\Delta_{\Gamma}$ the Laplace-Beltrami operator on $\Gamma$  and $\epsilon$ is   a positive constant representing  the  interface width.

The  origin of this kind of  partial differential equation resides in a more general Euler functional:
 \begin{equation}
 J_{\epsilon} (u)= \frac{1}{2}\int_{\Gamma} |\nabla u|^2 dA +  \int_{\Gamma} \frac{F(u)}{\epsilon^2}dA,
 \label{func0}
 \end{equation}
where
$$
 F(u)= \frac{1}{4}(u^2-1)^2.$$

This functional $J$ is nothing but the free  Helmholtz functional that then leads to  the known Allen-Cahn equation (\ref{AC1}).
We mention  other versions of the Helmholtz functional, such as logarithmic free  energy functional
$$
F(u) =\frac{\Theta}{2}[(1+u)\log(1+u) +(1-u)\log(1-u)]-\frac{\Theta_c}{2}u^2,
$$ where $\Theta$,$\Theta_c$  are constants) but they are not the focus of the current research.

We propose the following modified Allen-Cahn model (\ref{ModAC}):
\begin{equation}
u_t =\Delta_{\Gamma} u-\frac{1}{\epsilon^2}f_m(u), x\in \Gamma, t \in [0,T]
\label{ModAC}
\end{equation}
\begin{equation}
u|_{t=0}=u_0(x), x\in \Gamma
\end{equation}
where 
\begin{equation}
f_m(u)=u^3-u +b,
\label{f_mm}
\end{equation}
with
$b$ a real constant.

Its corresponding functional is a modified Helmholtz that  no longer has a double-well potential.
\begin{equation}
 J^m_{\epsilon} (u)= \frac{1}{2}\int_{\Gamma} |\nabla u|^2 dA +  \int_{\Gamma} \frac{F_m(u)}{\epsilon^2} dA,
 \label{func}
 \end{equation}
where 
$$
 F_m(u)= \frac{1}{4}(u^2-1)^2+b u $$
and $b$ is a real constant.

First we prove that this  modification does not change the functional property of being decreasing in  time. This  means that the total  energy  is a Lyapunov functional  for the solutions of the modified Allen-Cahn equation. 

\begin{Lemma}
The  energy functional  $J^m_\epsilon$  in (\ref{func})   is decreasing in time.
\end{Lemma}
Using Green's formula
$$\int_{\Gamma}  \nabla_{\Gamma} \xi \cdot \nabla_{\Gamma} \eta   dA = \int_{\partial \Gamma} \xi  \nabla_{\Gamma} \eta \cdot  \mu ds - \int_{\Gamma} \xi \Delta_{\Gamma} \eta dA.
$$
where $\mu$ is the  conormal on $\partial \Gamma$.
We then have by replacing $ \eta =u $ and $ \xi =u_t$ in the above formula and by means of integration by parts we get:

$$\frac{d J^m_{\epsilon}(u)}{dt} =\int_{\Gamma} \Big(\nabla_{\Gamma} u \cdot \nabla_{\Gamma} u_t + \frac{F_m'(u)}{\epsilon^2} u_t \Big) dA =$$

$$\int_{\Gamma} \Big(-\Delta_{\Gamma} u +\frac{F_m'(u)}{\epsilon^2}\Big)u_t dA =-\int_{\Gamma}u_t^2 dA  \leq 0.$$
This shows the decreasing behavior of the energy functional.
\subsection{One dimensional case analyis}

In this section we will   gain some intuition on the  quality  of the  underlying pattern  by understanding the 1D behavior. Assume for the one dimensional analysis that $\epsilon=1$. In the case when the spatial domain is an interval, we would look for stationary solutions, i.e $\frac{d u}{dt} =0$. Thus equation (\ref{ModAC})  becomes in 1D:
\begin{equation}
u"(x)+(u(x)-u^3(x)) -b =0
\label{1d}
\end{equation}

Here the potential is $F(u) = -\int f_m(u)du$.
with $f_m$ given in  (\ref{f_mm}).

Equation (\ref{1d}) can be solved analytically for $b=0$ , see \cite{nondegAC}. Its solution is given by
$u(x)=\tanh (x).$

For $b \neq 0$, multiplying (\ref{1d}) by $u'$:
\begin{equation}
u''u'-f_m(u)u' =0.
\end{equation}

Then
\begin{equation}
\frac{d}{dx}\Big(\frac{1}{2} \Big(\frac{du}{dx}\Big)^2+F(u)\Big)=0
\label{er}
\end{equation}

That means we can easily integrate (\ref{er}) as 
$$\frac{1}{2}(\frac{du}{dx})^2+F(u)=C$$

The solution is 
$$
\int_{u(a)}^{u(x)}\frac{du}{\sqrt{2(C-F(u))}} = x-a
$$

If we  use the exact expression of the potential $F(u)$,
the solution given by:
\begin{equation}
\int_{u(a)}^{u(x)}\frac{du}{\sqrt{2C-0.5(1-u^2)^2-2bu}}=x-a
\end{equation}
\begin{figure}[t]
\begin{center}
\includegraphics[width = 0.5\textwidth]{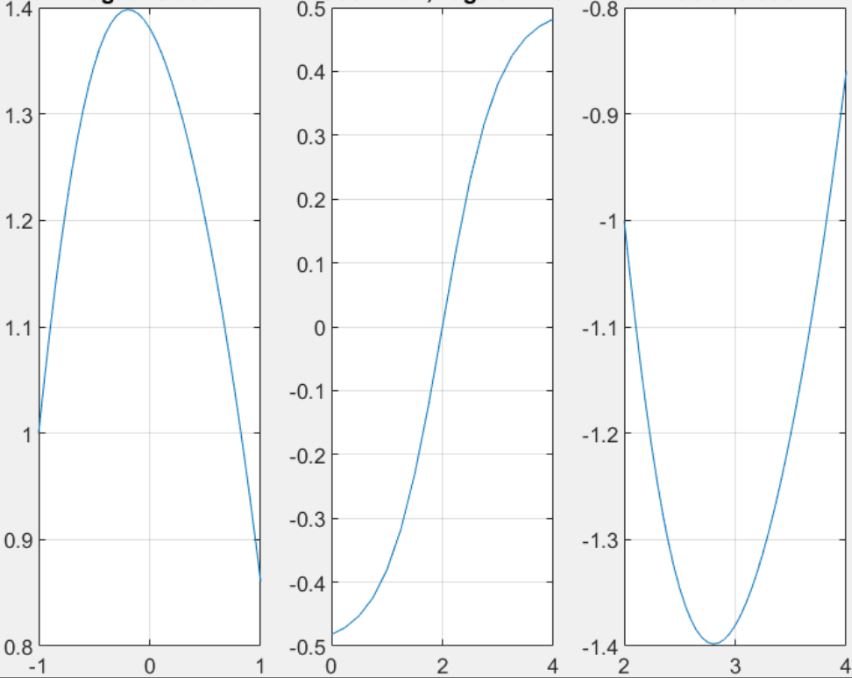}
\caption{Left. 1D solution of modified Allen-Cahn.$b = -1$. Middle. 1D solution to Allen-Cahn. $b = 0$. Right. 1D solution of modified Allen-Cahn. $b = 1.$}
\label{fig:1_dd}
\end{center}
\end{figure}
\begin{figure}[t]
\begin{center}
\includegraphics[width = 0.6\textwidth]{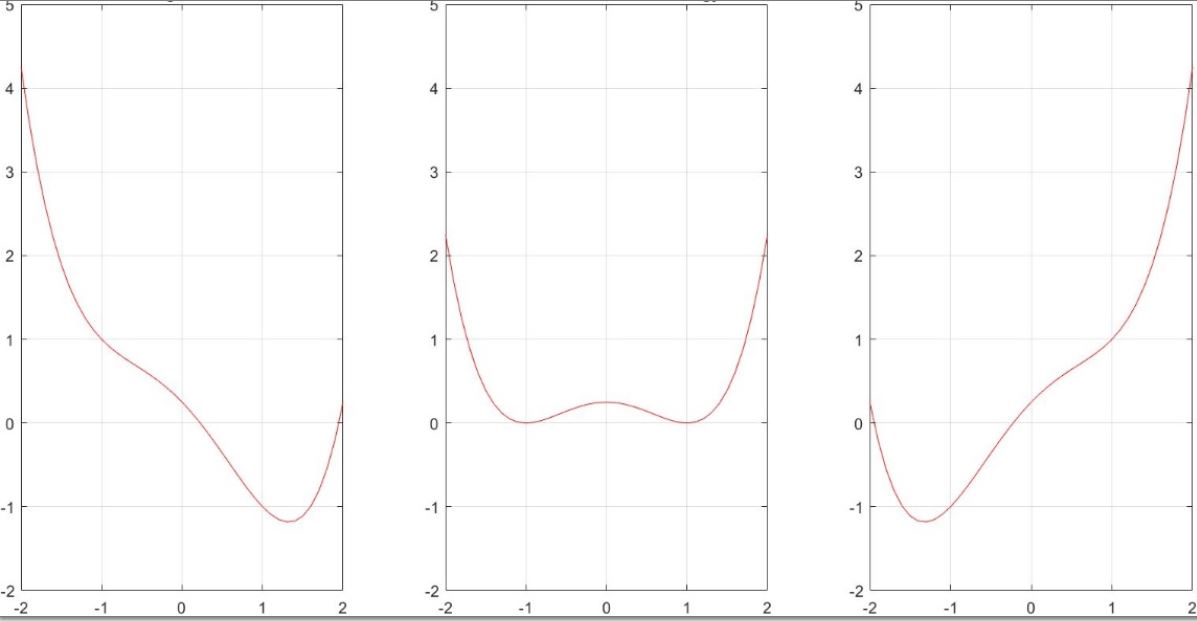}
\caption{Modified Allen-Cahn. Energy for the various reaction coefficients. Left. $b<0$. Middle. $b=0$. Right $b>0.$}
\label{fig:Graph_Helm}
\end{center}
\end{figure}

We used a solver  ode45 solver in Matlab to plot the numerical solution of the nonlinear ODE (for $b=1$ and $b=-1$):
$$
u"+ u(1-u^2)-b=0
$$
Even if we cannot assert  on patterns in the one dimensional case, in the graph Fig. \ref{fig:1_dd} one can observe the concave up form of the solution for positive coefficient  $b>0$, as well as  the concave down form of the solution for negative $b<0$.  In Fig. \ref{fig:Graph_Helm}  we plotted the corresponding Lyapunov energy.
As noticed in this simple one dimensional case, the solution is behaving in a clear relation with the sign of the reaction term. We will further illustrate for  more complex cases of various curvy  geometries, that the modified Allen-Cahn equation generates  spots/ inverted spots  in according to   negative/positive reaction coefficients.

We further  give details on the properties of  the  numerical scheme  for discretizing the proposed model.

\section{Numerical scheme}
{\bf Operator splitting based  scheme  for modified Allen-Cahn}

 The splitting of the operator emerges from the
 structure of the polynomial  structure in the reaction term that characterizes the specific structure of  the AC/modified AC equation.

Denote by $\Tilde{L}=L-\frac{b}{\epsilon^2} $ , where $L$ is the Laplace-Beltrami operator and by $B$ the nonlinear part ($Bu=\frac{u-u^3}{\epsilon^2}$).

According to Strang's splitting method ,  the numerical solution  to equation   in the time interval $[ t_n, t_{n+1}] $ can be written  as follows:

$$
U^{n+1} =(B^{\frac{\Delta t}{2}}\circ  \Tilde{L}^{\Delta t}\circ  B^{\frac{\Delta t}{2}}) U^n,
$$
We will write the above splitting operator into three steps:
\begin{equation}
 \tilde u =B(\tilde u), \tilde u^{n} =U^n, t\in [ t_n,t_{n+1}]
 \label{11}
\end{equation}

\begin{equation}
\bar u =\Delta_{\Gamma}\bar u -\frac{b}{\epsilon^2} , \bar u^{n} =\tilde u^{n+1},t\in [ t_n,t_{n+1}]
\label{22}
\end{equation}
\begin{equation}
\hat {u} =B( \hat {u}), \hat{u}^{n} =\bar u^{n+1}, t\in [ t_n,t_{n+1}]
\label{33}
\end{equation}
The numerical solution at $t=t_{n+1}$ is $ U^{n+1}=\hat{u}^{n+1}$.

The first and third  step  solve the same ODE, namely a Bernoulli equation, for which one can find   analytical solution.
\begin{equation}
\bar{u}^{n+1} =\frac{U^n}{\sqrt{e^{-\frac{2\Delta t}{\epsilon^2}}+(U^n)^2(1-e^{-\frac{2\Delta t}{\epsilon^2}})}}
\label{exAC}
\end{equation}
\textbf{Stability of the scheme}
\begin{Lemma}
For any time level $t=t_n$, the numerical solution $U^{n}$  given by (\ref{exAC}) for the first step  (\ref{11}) in the operator splitting is  unconditionally stable.
\end{Lemma}
\begin{Proof}
We have two possible cases.

Case A. $|U^n| \leq 1$, then
$$|\tilde u^{n+1}|=\frac{|U^{n}|}{\sqrt{ (U^n)^2 +(1-(U^n)^2) e^{-\frac{2\Delta t}{\epsilon ^2}}}} \leq \frac{U^n}{\sqrt{(U^n)^2}}=1$$
Case B. If $|U^n|>1$, then  again  one has
$$ |\tilde u ^{n+1}| \leq \frac{|U^n|}{ \sqrt {(1- e^{-\frac{2\Delta t}{\epsilon ^2}}) +  e^{-\frac{2\Delta t}{\epsilon ^2}}}} =|U^n|$$

Combining the two  cases, one can show

$$||\tilde  u ^{n+1}||_{\infty} \leq  \max\{|| U^n||_{\infty},1\},
$$
which completes the proof.
\end{Proof}

A similar Lemma can be formulated for the  third step  too.

 The second step of the splitting  involves discretization of the operator $\Tilde{L}$, i.e. the  Laplace Beltrami operator, for which we use cotan weight scheme for  triangulated meshes. The Laplace Beltrami  operator is  discretized  by $L =A^{-1}W$, where $A$ is the diagonal matrix of the Voronoi cells  areas around a vertex.
 \begin{figure}[h!]
\begin{center}
\includegraphics[width=30mm]{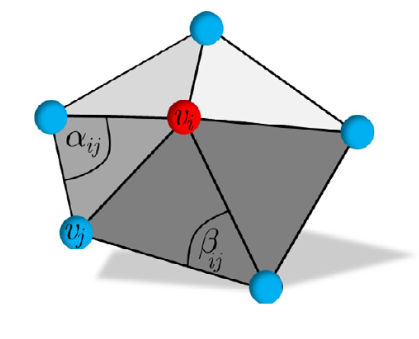}
\caption{Weights in Discrete Laplacian}
\label{fig:20_DLB}
\end{center}
\end{figure}
$$ W_{ij} =\begin{cases} 
     \sum_{v_j\in N_i} w_{ij} ,i =j \\
      -w_{ij}, i\neq j, v_j \in N_i \\
       0 ,\quad  \mbox{otherwise}
   \end{cases}
   $$
and the weights $w_{ij} = \cot  \alpha _{ij}+ \cot \beta_{ij}$  where $\alpha_{ij}$ and $\beta_{ij}$ are the angles opposite to the edge  as appearing in Fig \ref{fig:20_DLB}. If a small  step time parameter is used, step two of the operator splitting scheme  is  stable.
\begin{Remark}
The operator splitting scheme (\ref{11}), (\ref{22}),(\ref{33}) is stable only under small time-step restriction. To obtain a fully unconditionally  operator  splitting scheme, one might further use a full FEM  or a Cranck-Nicolson scheme for discretizing the part involving the Laplace-Beltrami operator.
\end{Remark}
\section{Numerical examples}
\subsection{Global Patterns with modified Allen-Cahn model}
In this section  we will show  numerical  experiments   of phase separation as a result of    applying   Allen-Cahn model as well as  modified  A-C model on various surfaces to show  the generated patterns.
While with Allen-Cahn PDE, the only  pattern that can be obtained is a stripe based  shape, the  modified new model  allows generating  various patterns in form of spots/inverted spots  on the surfaces and an  operator-splitting scheme is used to easily  implement it.

\begin{figure}[!htb]
\begin{center}
\includegraphics[width=110mm]{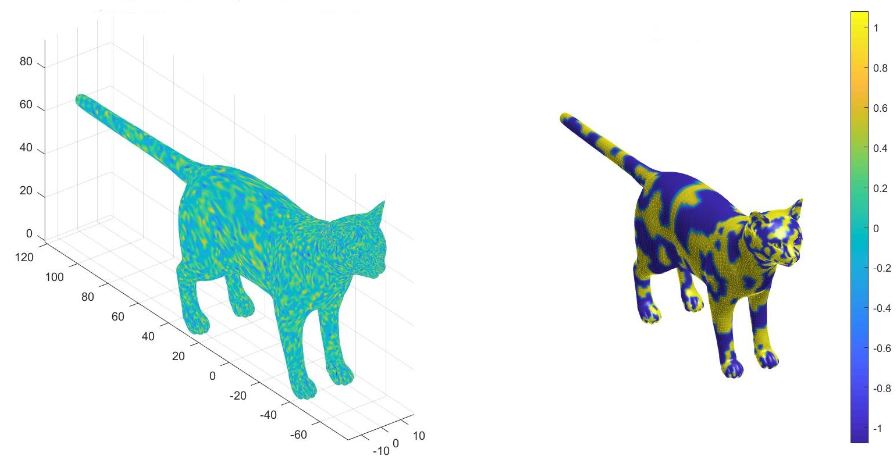}
\caption{ Allen-Cahn on an animal  surface. Left: Random  initial  data. Right: Stripes by Allen-Cahn. Reaction term:  $b = 0. $}
\label{fig:2}
\end{center}
\end{figure}

\begin{figure}[!htb]
\begin{center}
\includegraphics[width=120mm]{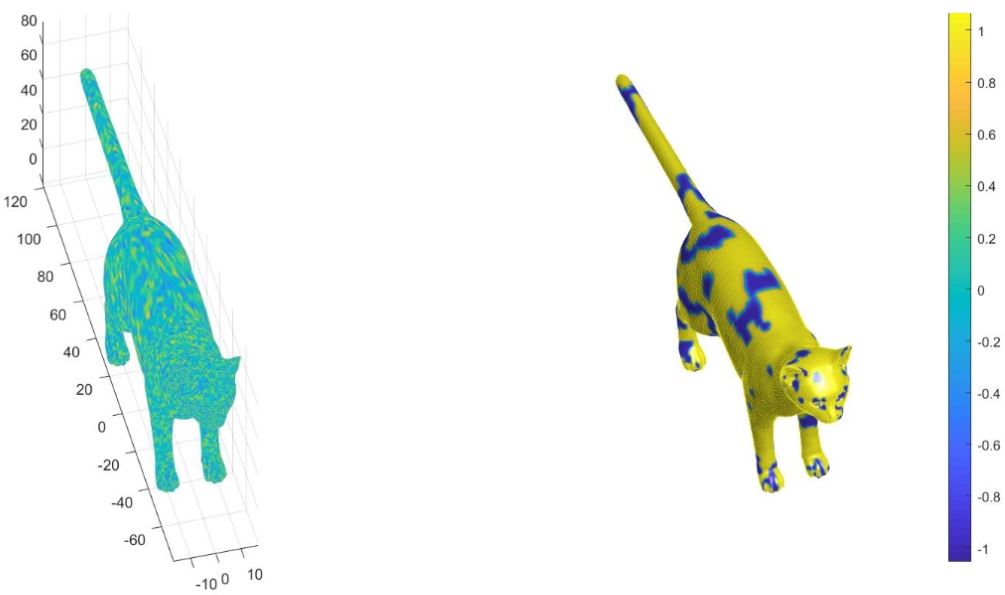}
\caption{Left: Random  initial  data. Right: Inverted spots. Reaction term  : $b = 0.2 $.}
\label{fig:3}
\end{center}
\end{figure}

\begin{figure}[!htb]
\begin{center}
\includegraphics[width=120mm]{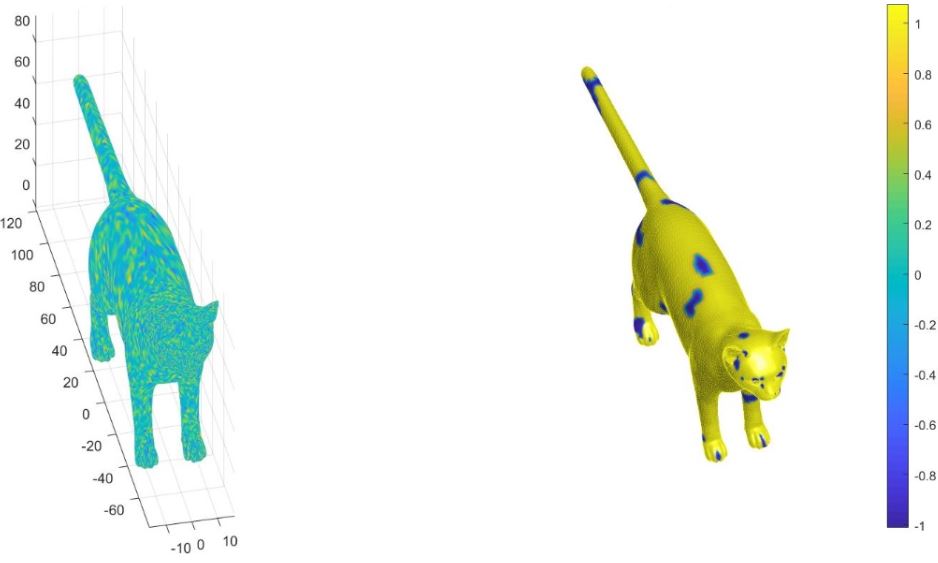}
\caption{Left: Random  initial  data. Right: Inverted spots. Reaction term  : $b = 0.3 $.}
\label{fig:3b}
\end{center}
\end{figure}

\begin{figure}[!htb]
\begin{center}
\includegraphics[width=120mm]{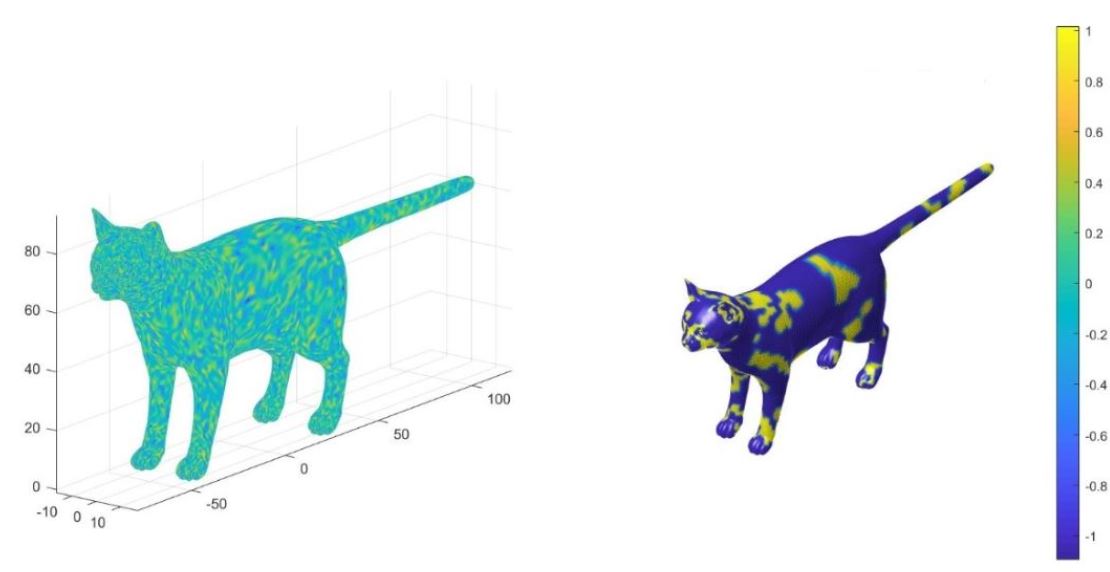}
\caption{Left: Random  initial  data. Right: Spots.  Reaction term.  $b = - 0.2$.} 
\label{fig:4}
\end{center}
\end{figure}

\begin{figure}[!htb]
\begin{center}
\includegraphics[width=120mm]{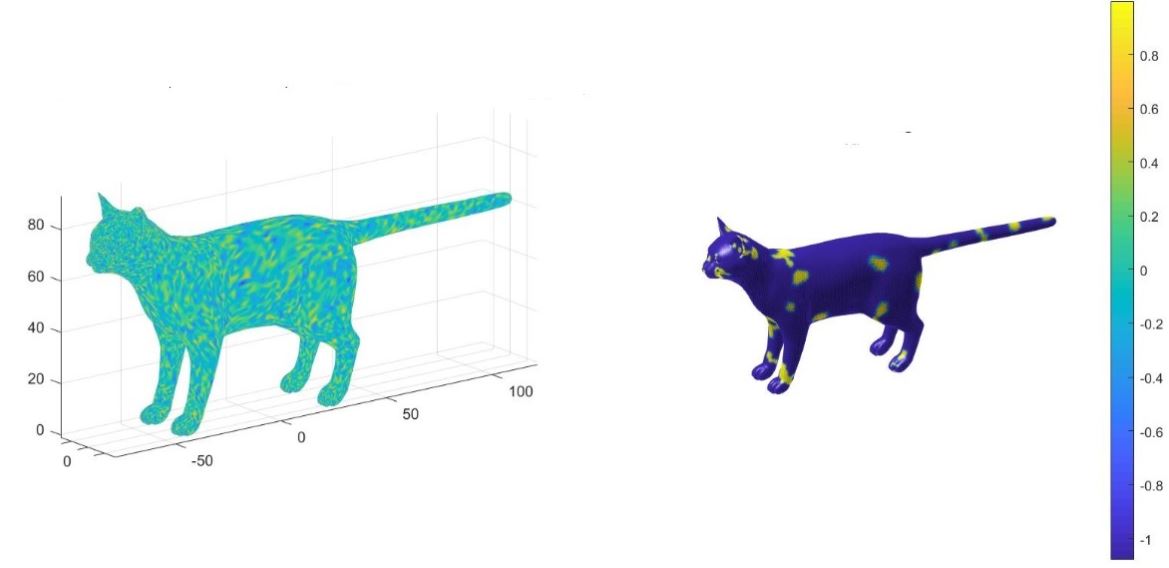}
\caption{Left: Random  initial  data. Right: Spots.  Reaction term.  $b = - 0.3$.} 
\label{fig:4b}
\end{center}
\end{figure}

\begin{figure}[!htb]
\begin{center}
\includegraphics[width=120mm]{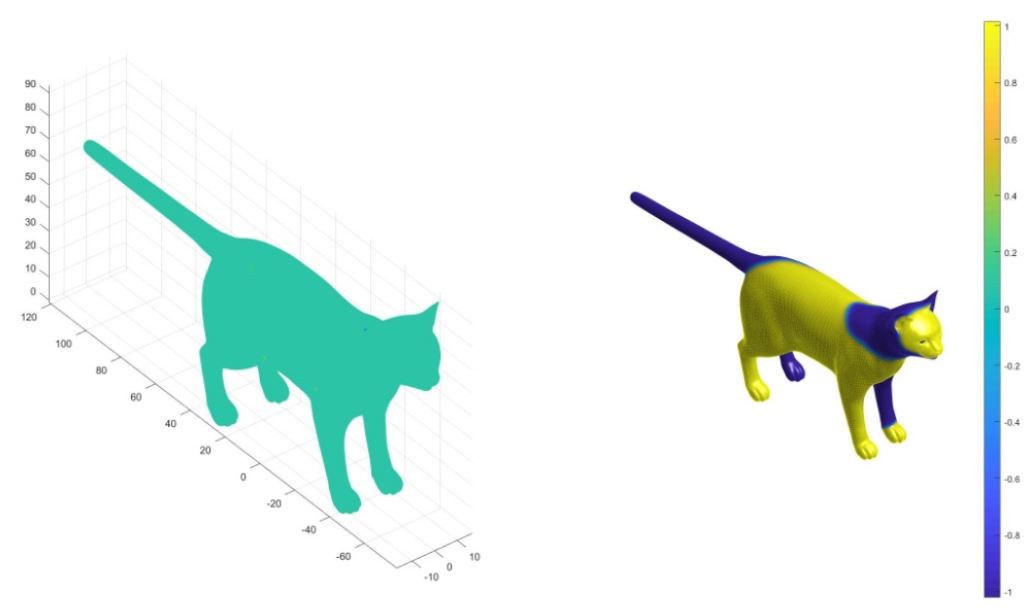}
\caption{  Left: Initial sparse random data. Right: The stripes created with Allen-Cahn model.}
\label{fig:6}
\end{center}
\end{figure}

The characteristics of the  resulting patterns obtained by applying the modified  Allen-Cahn  (\ref{ModAC}) model are  determined   according  to the values of the constant $b$.
The following types of patterns can be generated (see Fig.\ref{fig:2}-Fig.\ref{fig:6}):
for $b>0$  inverted spots, for $b<0$ spots, and for
$b=0$, stripes,  i.e. regular Allen-Cahn.  
The generated pattern can be categorized, as seen in the numerical experiments on various surfaces, as follows:
 when $0<b<0.3$, one obtains inverted spots,
 when $-0.3<b<0$ one gets spots, while 
 for values $|b|>0.3$, the diffusion term takes over the reaction one, and one gets a trivial  constant solution as time progresses.
 
\vspace{4cm}
%\vspace{1cm}
\clearpage

\begin{figure}[!htb]
\begin{center}
\includegraphics[width=120mm]{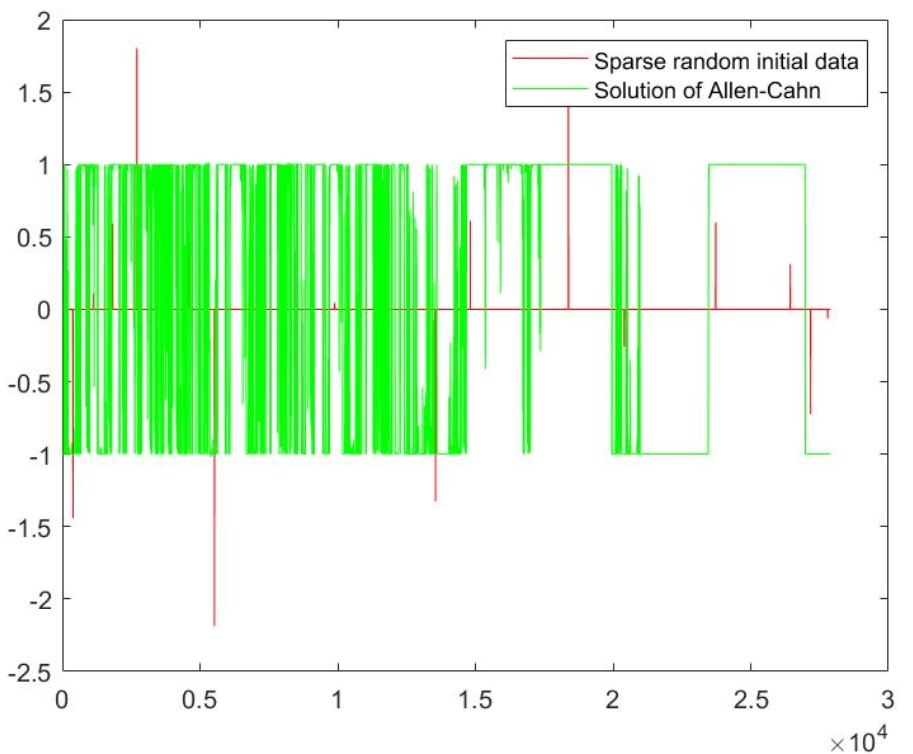} 
\caption{Graph of initial  data and of numerical solution to  Allen-Cahn  model.}
\label{fig:6b}
\end{center}
\end{figure}
The example in Fig.\ref{fig:6} illustrates  the stripes created by Allen-Cahn    model, given a sparse random initial data on the surface. The solution  exhibits the phase separation property, see Fig.\ref{fig:6b} i.e. the solution tends to  $\pm1$. An interface of width $\epsilon$  is  created  between the two phases.

\begin{figure}[!htb]
\begin{center}
\includegraphics[width=100mm]{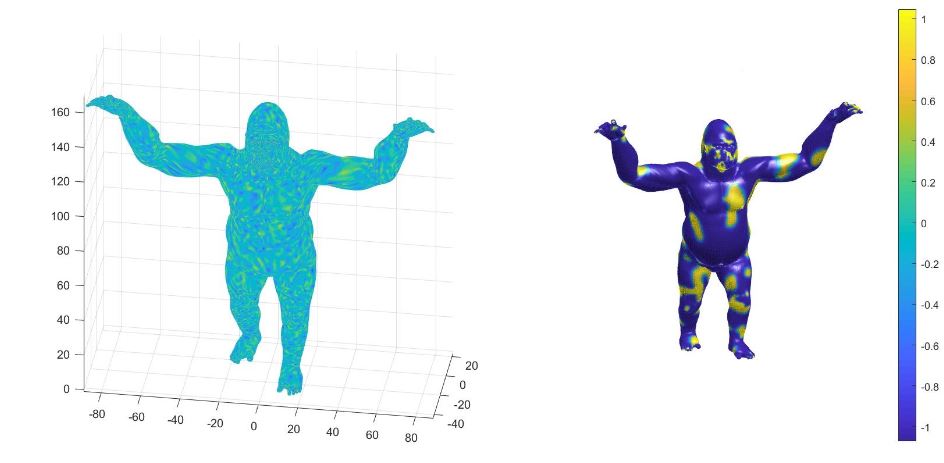}
\includegraphics[width=100mm]{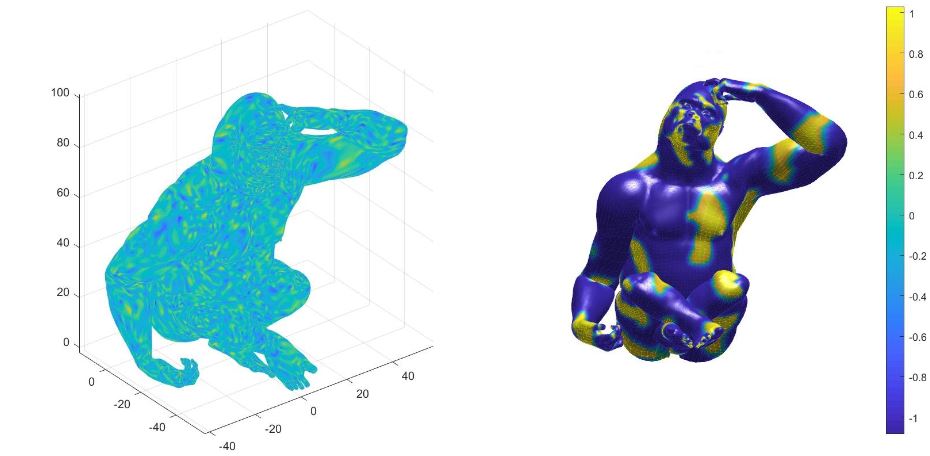}
\includegraphics[width=100mm]{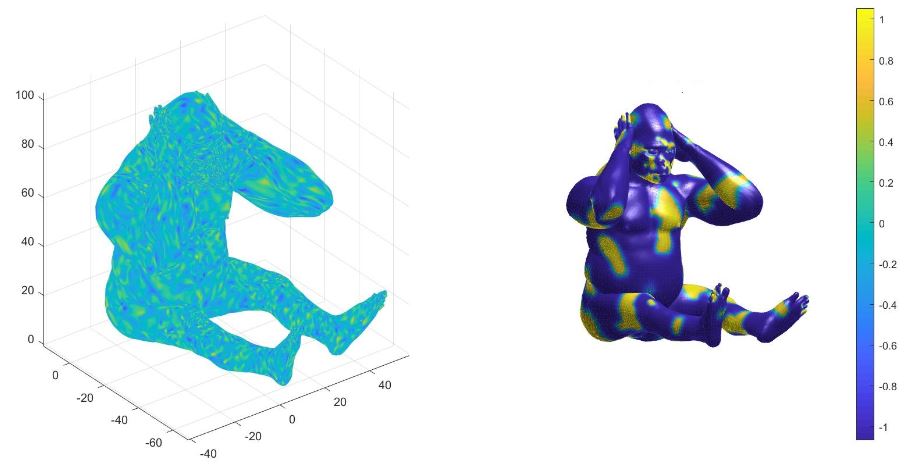}
\caption{Various almost isometric shapes with generated global spot  patterns in corresponding  regions on the surfaces. Giving the same initial  random data defined on  surface, the modified A-Cahn model  leads to    spots in the  three isometric surfaces located in  corresponding regions.}
\label{fig:70}
\end{center}
\end{figure}

\begin{figure}[!htb]
\begin{center}
\includegraphics[width=80mm]{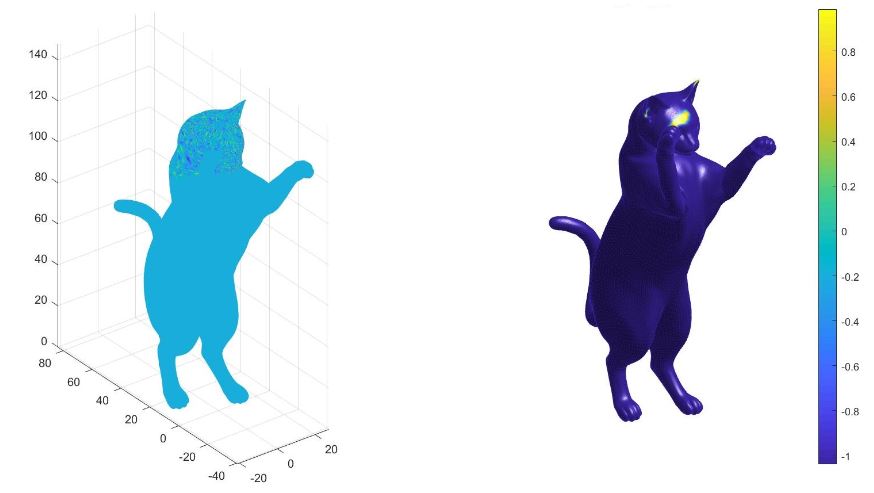}
\includegraphics[width=85mm]{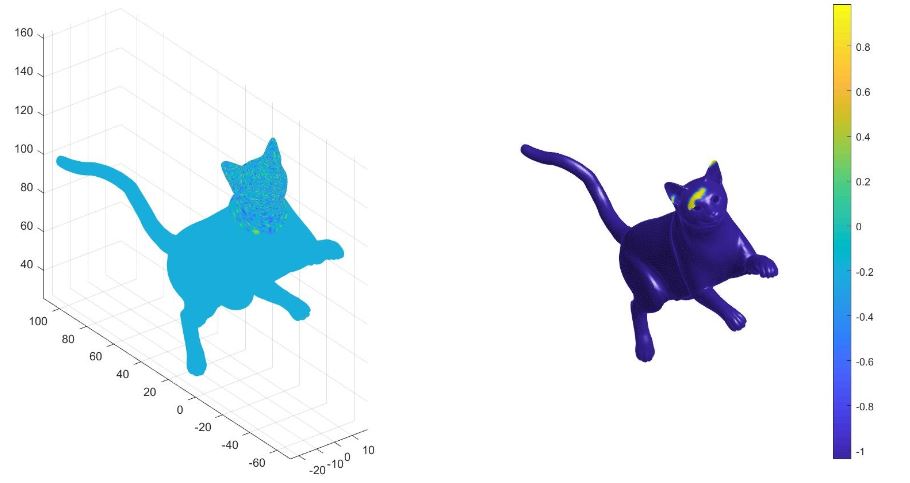}
\includegraphics[width=90mm]{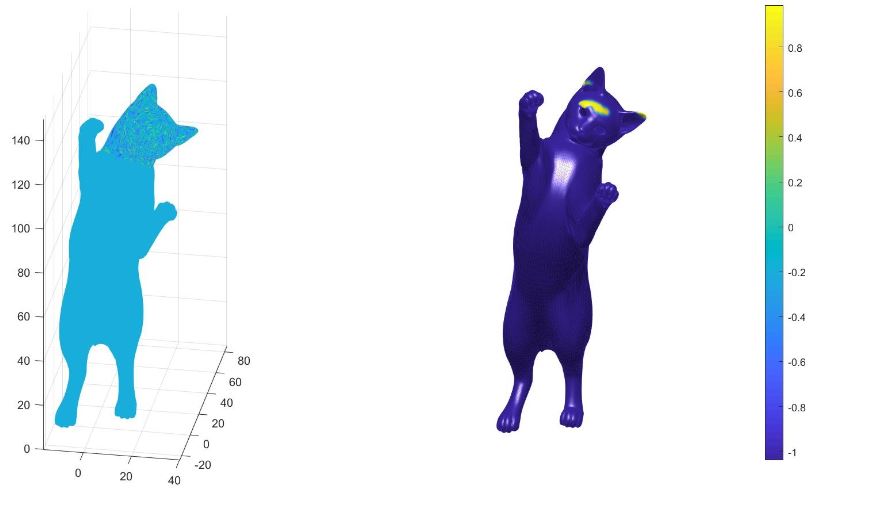}
\caption{Various almost isometric shapes with generated local  spot  patterns in corresponding  regions on the surfaces. Giving the same initial  random data defined on compact set on the surfaces, the modified A-Cahn  model leads to   localized  spots residing in the same compact set, in the  three isometric surfaces in  corresponding regions.  $b = -0.08.$ }
\label{fig:88}
\end{center}
\end{figure}

The examples  in Fig.\ref{fig:70}, Fig.\ref{fig:88} show  examples of generated spots for almost invariant surfaces.   Giving the same initial  random data,  the modified Allen-Cahn model leads to    spots in the  three isometric surfaces located in  corresponding regions. This  invariance property resides in  the fact that the modified Allen-Cahn is invariant to isometries due to the  Laplace-Beltrami operator which  is  invariant to isometries,  and moreover the nonlinear part of the reaction term of the governing equation  depends only on $u$, and not on the local  coordinates on the surface. 
\newpage
\subsection{Localized patterns with modified Allen-Cahn model}
\begin{figure}[!htb]
\begin{center}
\includegraphics[width=100mm]{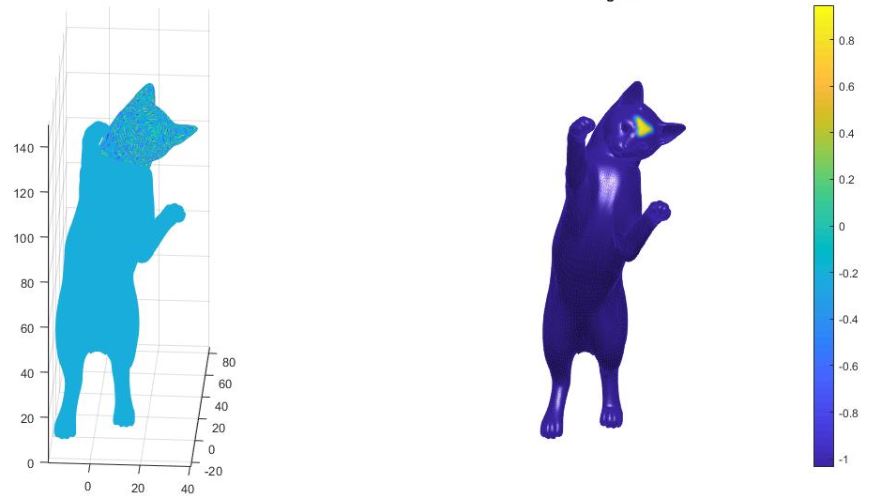}
\includegraphics[width=100mm]{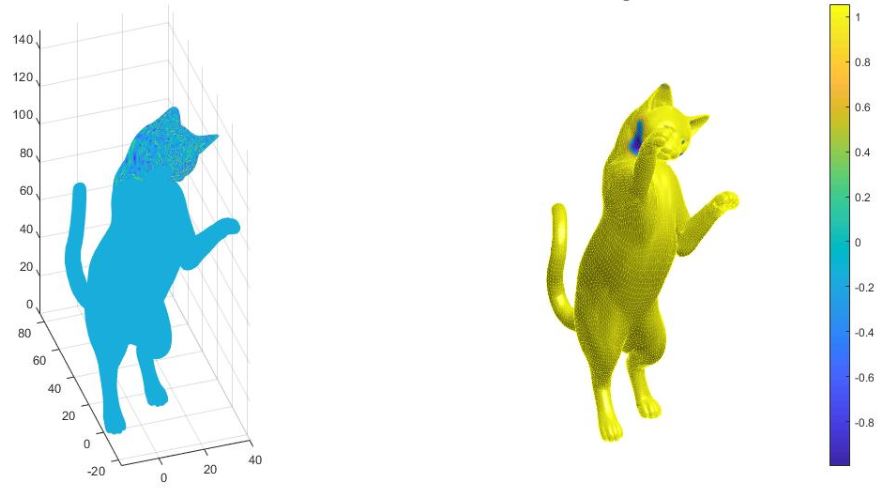}
\caption{Left top and bottom. Localized random initial data.Top right:  Generated spots  remain localized, residing  in the area where the initial data was defined. $b$ = -0.08. 1200 iterations. $dt $= $0.9$. Bottom right. Inverted spots. $b = 0.08$.}
\label{fig:8}
\end{center}
\end{figure}

\begin{figure}[h!]
\begin{center}
\includegraphics[width=120mm]{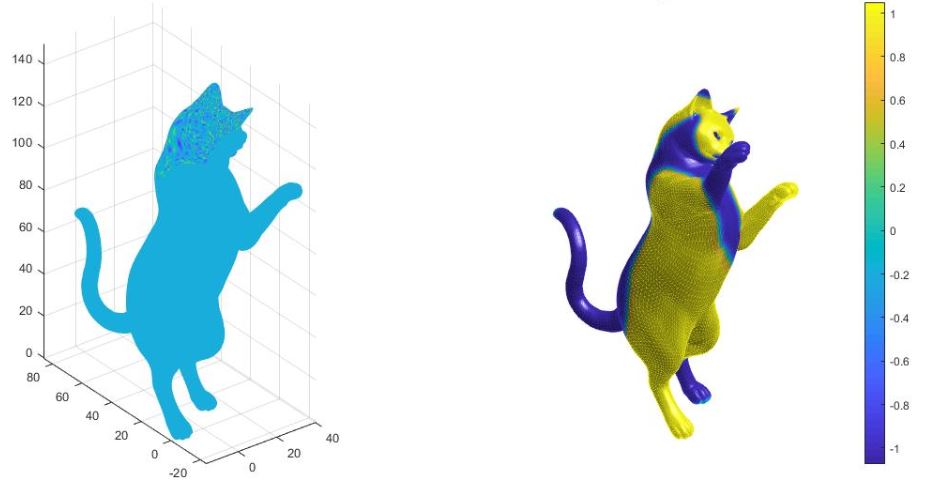}
\includegraphics[width=120mm]{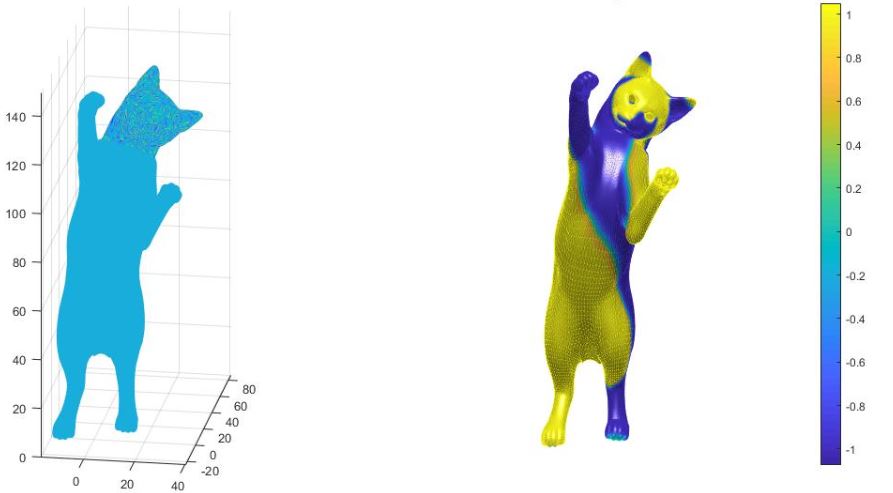}
\caption{Left top and bottom. Localized  random  initial data.  Right. Generated stripes (Allen-Cahn, $b = 0$), see two views on the pattern on the surface. Generated pattern is not local, moulding on the whole surface.
 }
\label{fig:9}
\end{center}
\end{figure}

The locality of the modified model is an important property and is illustrated in examples  below.
While in the case of Allen-Cahn it fails see Fig.  \ref{fig:9}, this locality  is satisfied only by the modified Allen-Cahn model, see Fig. \ref{fig:8}. Given an initial data in a compact domain, by means of the modified A-Cahn model, the generated pattern will be formated only in this domain.

\clearpage
\clearpage

\section{Conclusions} 
We propose a slight modification of the known reaction-diffusion Allen-Cahn model on surfaces.
  Unlike with the Allen-Cahn model,  generating only stripes,  the  modified model  can be used to generate patterns in form of spots/inverted spots. The dependence of the   kind of pattern on the new  introduced reaction  term is investigated. Furthermore, a simple and efficient operator based splitting scheme is employed to  discretize the equation. Numerical examples show the solution of the corresponding PDE under varying  initial conditions to illustrate   various patterns and the underlying local or global  generated patterns.  In future work  we will include  extensions of the  modified  Allen-Cahn  flow on deforming geometries and exploring   possible  applications of such  patterns   for shape analysis/correspondence.
\bibliography{ACahnModifModel}
\bibliographystyle{unsrtnat}  
\end{document}